\documentclass[12pt]{amsart}
\input diagrams
\diagramstyle[PostScript=dvips]

\newcommand{\cal}{\mathcal}
\renewcommand{\Bbb}{\mathbb}

\newcommand{\ch}{\operatorname{ch}}
\newcommand{\td}{\operatorname{td}}

\newcommand{\DD}{{\cal D}}

\newcommand{\MM}{{\cal M}}
\newcommand{\mbar}{\overline{\MM}}
\newcommand{\mgn}{\overline{\MM}_{g,n}}
\newcommand{\mgnr}{\mgn^{1/r}}

\newcommand{\cv}{c^{1/r}}            

\newcommand{\EE}{{\cal E}}

\newcommand{\TT}{{\cal T}}
\newcommand{\XX}{{\cal X}}

\newcommand{\Bl}{\operatorname{Bl}}

\newcommand{\hra}{\hookrightarrow}
\newcommand{\lan}{\langle}
\newcommand{\ran}{\rangle}

\newcommand{\CC}{{\cal C}}

\renewcommand{\P}{{\Bbb P}}

\newcommand{\si}{\sigma}

\newcommand{\We}{\Lambda}
\renewcommand{\ker}{\operatorname{ker}}
\newcommand{\im}{\operatorname{im}}

\numberwithin{equation}{section}

\newtheorem{thm}{Theorem}[section]
\newtheorem{prop}[thm]{Proposition}
\newtheorem{lem}[thm]{Lemma}
\newtheorem{cor}[thm]{Corollary}
\newenvironment{rem}{\vspace{3mm}\noindent
{\bf Remark.}}{\vspace{3mm}}
\newenvironment{rems}{\vspace{3mm}
\noindent {\bf Remarks.}}{\vspace{3mm}}

\theoremstyle{definition}

\newcommand{\Pf}{\noindent {\it Proof}}
\newcommand{\id}{\operatorname{id}}

\newcommand{\ov}{\overline}

\newcommand{\rk}{\operatorname{rk}}
\newcommand{\ra}{\rightarrow}
\newcommand{\isomarrow}{\rTo^{\sim}} 

     \renewcommand{\SS}{{\cal S}}

\newcommand{\OO}{{\cal O}}

\newcommand{\Hom}{\operatorname{Hom}}
\newcommand{\Ext}{\operatorname{Ext}}

\renewcommand{\a}{\alpha}
\renewcommand{\b}{\beta}
\newcommand{\om}{\omega}

\newcommand{\la}{\lambda}
     
\newcommand{\C}{{\Bbb C}}

\newcommand{\Z}{{\Bbb Z}}
\newcommand{\Q}{{\Bbb Q}}

\newcommand{\Ga}{\Gamma}

\newcommand{\wt}{\widetilde}

\newcommand{\gr}{\operatorname{gr}}
\newcommand{\sub}{\subset}
\newcommand{\ed}{\qed\vspace{3mm}}
\newcommand{\A}{{\Bbb A}}

\begin{document}

\title[Witten's top Chern class]
{Algebraic construction of Witten's top Chern class}

\author
[A.\,Polishchuk]{Alexander Polishchuk}
\address
{Department of Mathematics, Boston
University, Boston, MA 02215} 
\email{apolish@math.bu.edu}
\thanks{Work of the first author was partially supported by NSF grant
DMS-0070967.}

\author
[A.\,Vaintrob]{Arkady Vaintrob}
\address
{Department of Mathematics, University of Oregon, Eugene, OR 97403}
\email{vaintrob@math.uoregon.edu}

\date{October, 2000}

\begin{abstract} 

Applying a modification of MacPherson's graph construction
to the case of periodic complexes, we give an algebraic construction of
Witten's ``top Chern class'' on the moduli space of algebraic curves with
higher spin structures. We show that it satisfies most of the axioms 
for the spin virtual class and also the so-called descent axiom. 
\end{abstract}
\maketitle

\section{Introduction}

A famous conjecture of Witten~\cite{W3} proved by Kontsevich
in~\cite{Ko} connects the generating
function of certain intersection numbers on the moduli space
$\mgn$ of $n$-pointed stable curves of genus $g$ and a special
solution of the Korteveg-de-Vries hierarchy. 
In an attempt to generalize this conjecture to the so-called
Gelfand-Dickey hierarchies, Witten in~\cite{W1} 
introduced a certain physical model of two-dimensional quantum gravity
coupled to matter. He interpreted the free energy of this model
in terms of intersection theory on certain ramified covers of 
 $\mgn$ obtained by compactifying the moduli
space of Riemann surfaces with $r$-th roots of the canonical bundle
($r$-spin structures) 
and conjectured that it gives a special solution of
the $r$-th Gelfand-Dickey hierarchy. Unlike the original conjecture,
however, the  generalized conjecture
was not precisely formulated at that time.
 First, the appropriate compactification of the space  
of curves with $r$-spin structure had not yet been constructed.
 In~\cite{W2},  
Witten described the desired types of behavior of a spin structure 
near the singular points and Jarvis in~\cite{J2,J}  gave a detailed
algebro-geometric construction of the   
corresponding moduli space $\mgnr$. Second, the intersection numbers
in question involved not only the Chern classes of the tautological
line bundles, but also the so-called top Chern class $\cv$ of
the push-forward (in the derived category sense)
of the universal $r$-th root bundle.
Since this push-forward is in general not a sheaf, but only a complex
of sheaves, the naive definition of the class $\cv$ does not work.
This resembles the situation in the definition of the Gromov-Witten
invariants of a variety $V$, where one has to integrate certain
cohomology classes against the fundamental class of the moduli space
of stable maps from Riemann surfaces to $V$. Since the space of
stable maps is  not smooth in general, the non-existent ``honest''
fundamental class has to be replaced by a so-called ``virtual''
fundamental class. This virtual class admits two different constructions,  
analytic and algebraic, both of which satisfy a set of axioms 
formulated by Behrend and Manin in~\cite{BM}.

Similar axioms for Witten's {\em virtual top Chern class\/}
 $\cv$ on the moduli space $\mgnr$
were formulated in~\cite{JKV}; a class satisfying these
axioms is called a {\em spin virtual class\/}.
It was proven there that any spin virtual class gives
rise to a cohomological field theory (in the sense of
 Kontsevich-Manin~\cite{KM}). Existence of such class (as well as
the generalized Witten's conjecture) was also proven in~\cite{JKV}
in two special cases:  $g=0$ (and arbitrary $r$) and $r=2$ (and arbitrary
 $g$). 

In~\cite{W2}, Witten sketched an index-theoretic construction
of the virtual top Chern  class. However, it is not clear how to apply 
his construction  on the boundary strata of $\mgnr$ and whether the
resulting class would have the desired properties.
\medskip

In this paper we give an algebraic construction of a Chow cohomology class on
$\mgnr$ that we propose as a candidate for the role of the spin virtual class.
Our class satisfies all axioms  of the spin virtual class
formulated in~\cite{JKV} except, perhaps,  some
vanishing conditions which we verified only in several special cases.
For the cases considered in~\cite{JKV} our construction gives the same
answer.
\medskip

The paper is organized as follows.
In Section 2 we describe a modification of MacPherson's 
graph construction for the case of 2-periodic (unbounded) complexes.
In Section 3 we use it to define
an analogue of the localized top Chern class for an orthogonal bundle $E$ with
an isotropic section $s$. The corresponding periodic complex 
is the spinor bundle of $E$ where the differential comes from
an action of $s$. In Section 4 we apply this construction to orthogonal
bundles related to families of $r$-spin structures and in Section 5 we 
study properties of this virtual top Chern class. 

Our construction can be viewed as an algebraic counterpart of Witten's
analytic construction. It would be interesting to fill in the details
of  Witten's construction and to check whether the resulting classes
coincide. (In the case of the virtual fundamental class, it is known
that algebraic and analytic constructions give the same result).
\medskip

\noindent
{\em Acknowledgment.\/}
We are grateful to the Institute des Hautes \'Etudes Scientifiques
for its hospitality during the final stages of work on this paper.
Also we would like to thank Takashi Kimura for reading the first
draft of the paper and for making valuable comments.
\medskip

\noindent
{\em Conventions.\/} 
All our schemes are assumed to be of finite type
over $\C$. By a variety we mean a reduced and irreducible scheme.
When we consider the zero locus of a section of a vector bundle,
it is assumed to be  equipped with the reduced scheme structure.

\section{Periodic graph-construction}

\subsection{Bivariant intersection theory}

We refer to Fulton's book \cite[Ch.\ 17]{F} for the main definitions  
concerning bivariant intersection theory. 
The following lemma combined with the deformation to the normal cone
(see \cite[Ch.\ 5]{F}) shows that in order to prove some universal
identities between bivariant classes involving a closed embedding, it
is sufficient to prove these identities in the case of an embedding
of a zero section into a normal cone.

\begin{lem}\label{deform}
Let $U\subset\P^1$ be an open subset, containing the point
$\infty$, $p:\XX\ra U$ be a flat morphism such that there
is an isomorphism $p^{-1}(U\setminus\{\infty\})
\simeq X\times (U\setminus\{\infty\})$ of $U\setminus\{\infty\}$-schemes
where $X$ is a scheme. For every $t\in U$ let us denote by
$i_t:\XX_t\hra\XX$ the fiber of $\XX$ over $t$ (so that $\XX_t=X$
for $t\neq\infty$).
Let $Z\times U\ra \XX$ be a family of closed embeddings $k_t:Z\hra\XX_t$
parametrized by $U$.
Assume that we have a bivariant class $c\in A^n(Z\times U\ra\XX)$
such that $i_{\infty}^*c=0$ in $A^n(Z\ra\XX_{\infty})$. Then
$i_t^*c=0$ in $A^n(Z\ra\XX_t)$ for all $t\in U$.
\end{lem}

\Pf . It is sufficient to check that for every $t\in U\setminus\{\infty\}$
and every class $\a\in A_p(\XX_t)$ one has
$c(\a)=0$ in $A_{p-n}(Z)$. Let
$$j:X\times (U\setminus\{\infty\})\simeq p^{-1}(U\setminus\{\infty\})\ra\XX$$
be the natural open embedding. 
We can choose a class $\wt{\a}\in A_{p+1}(\XX)$ such that
$j^*\wt{\a}=p^*\a$ in $A_{p+1}(X\times (U\setminus\{\infty\}))$.
Since bivariant operations commute with specialization we have
$$c(\a)=i_t^!c(\wt{\a}).$$
By example 2.6.6 of \cite{F}, the Gysin homomorphism
$A_*(Z\times U)\ra A_*(Z\times \{t\})$ doesn't depend on $t$.
Hence, we have
$$i_t^!c(\wt{\a})=i_{\infty}^!c(\wt{\a})=c(i_{\infty}^!\wt{\a})=0$$
since $i_{\infty}^*c=0$.
\ed

\subsection{Localized Chern character for $2$-periodic complexes}

Recall that the graph-construction (see \cite{BFM} or \cite[Ch.\ 18]{F})
associates with a finite complex of vector bundles $V^{\bullet}$ on $X$
which is exact off a closed subset $Z\subset X$, the
localized Chern 
character
$$\ch_Z^X(V^{\bullet})\in A(Z\ra X)_{\Q}.$$
Notice that every finite complex $V^{\bullet}$
defines an infinite $2$-periodic complex:
$$\ldots \ra V^+\rTo^{d^+} V^-\rTo^{d^-} V^+\ra\ldots$$
where 
$$V^+:=\bigoplus_{n\in\Z} V^{2n}, \quad V^-:=\bigoplus_{n\in\Z} V^{2n+1},$$ 
and the maps $d^+$ and $d^-$ are induced by the differentials in $V^{\bullet}$.
We observe that the localized Chern character
can be defined starting from this $2$-periodic complex
so one can generalize it to the case of infinite complexes.
The only subtlety in the case of infinite complexes is that one
must be more careful about the condition of exactness off $Z$. 

Namely, we say that a $2$-periodic
complex $V^{\bullet}$ is {\em strictly exact off $Z$\/}
if $V^{\bullet}$ is exact off $Z$ and the images of $d^+$ and $d^-$
are subbundles. If $X$ is a variety this is equivalent to the condition that 
$V^{\bullet}|_x$ is exact for every 
point $x\in X\setminus Z$. 
\medskip

We can consider a $2$-periodic complex as a $\Z/2\Z$-graded vector bundle
$V=V^+\oplus V^-$ on $X$ 
equipped with homomorphisms $d^+:V^+\ra V^{-}$ and 
$d^-:V^-\ra V^+$ such that $d^+d^-=d^-d^+=0$. Let $v^+=\rk V^+$ and
$v^-=\rk V^-$. 

Denote by $G^+$ (resp.\ $G^-$) the Grassmannian
bundle of $v^+$-dimensional (resp.\ $v^-$-dimensional)
planes in $V$. Let $G=G^+\times_X G^-$ and let $\xi^+$ and 
$ \xi^-$ be the pull-backs of the tautological bundles on $G^+$ and $G^-$
under the natural projection $G\ra G^\pm$ respectively. 
Denote by $\xi$ the virtual bundle $\xi=\xi^+-\xi^-$ on $G$.

We have a natural embedding
$$\varphi:X\times\A^1\rInto G\times\A^1:
(x,\la)\mapsto(\gr(\la d^+(x)),\gr(\la d^-(x)),\la)$$
where by $\gr(f)$ we denote the graph of a linear map $f$.
There is also a second
 natural embedding
$$\psi:(X\setminus Z)\times\P^1\ra G\times\P^1$$
that sends a point $(x,(\la_0:\la_1))$ to 
$(H^+,H^-,(\la_0:\la_1))$ where subspaces
$$ H^+ \subset V^+\oplus\ker d^- \ \text{ \ and \ } H^- \subset \ker d^+ 
\oplus V^- $$   
are defined by 
$$H^+=\{(x_0,x_1)\,|\,\la_0x_1=\la_1 d^+(x_0)\}
 \text{\ and \ } H^-=\{(x_0,x_1)\,|\,\la_0x_0=\la_1 d^-(x_1)\}.
$$
The requirement that the complex $V^{\bullet}$ is strictly exact off $Z$
is essential to
define these embeddings scheme-theoretically.
\medskip

The localized Chern character
$$\ch_Z^X(V^{\bullet})\in A^*(Z\ra X)_{\Q}$$ 
for a strictly exact periodic complex can now be defined
using the data $(\phi,\psi,\xi)$ in the same way as in \cite[Ch.\ 18]{F}.
The key fact needed for the construction is the following property.

Consider the embedding
$$\psi^0:(X\setminus Z)\times\A^1\rInto G\times\A^1$$
induced by $\psi$. Then the pull-back of the virtual bundle
$\xi$ by $\psi^0$ is trivial.
\medskip

In the case $Z=X$ this construction gives the usual graded Chern character
of $V^{\bullet}$:
$$\ch_X^X(V^{\bullet})=\ch(V^{\bullet}):=\ch(V^+)-\ch(V^-).$$

\subsection{Properties of $\ch_Z^X(V^{\bullet})$}

Here we will establish some basic properties
of the graph-construction that will be used later.
The first property is the compatibility of our construction 
with the standard one in the case of bounded complexes.

\begin{prop} Let $(C^{\bullet},d)$ be a bounded complex
of vector bundles on $X$ strictly exact off $Z$ and
$C^{\bullet}_{(2)}$ be the corresponding $2$-periodic complex given by 
$$C_{(2)}^+=\bigoplus_{i\equiv 0 (\mod 2)}C^i, \qquad
C_{(2)}^-=\bigoplus_{i\equiv 1 (\mod 2)}C^i$$
with the differentials $d^{\pm}$ induced by $d$. Then one has
$$\ch_Z^X(C^{\bullet})=\ch_Z^X(C^{\bullet}_{(2)}).$$
\end{prop}

\medskip

{}From now on we will consider only $2$-periodic complexes. When we say
that such a complex is exact off $Z$ we will always mean strict
exactness (i.e.\ that $\im(d^+)$ and $\im(d^-)$ are subbundles).
\medskip

As in the case of bounded complexes one can prove the following properties
of the periodic graph-construction.

\begin{prop}\label{properties}
Let $V^{\bullet}$ be a $2$-periodic complex of vector bundles
on $X$ exact off $Z$. 

\begin{enumerate}
\item[(i)] Let $i:Z\hra Z'$ be a closed embedding where $Z'\subset X$
is another closed subset of $X$. Then
$$i_*\ch_Z^X(V^{\bullet})=\ch_{Z'}^X(V^{\bullet}).$$

\item[(ii)]
 Let $U\subset X$ be an open subset such that $Z\subset U$
and $[j]$ the canonical orientation class of the embedding $j:U\rInto X$.
Then
$$\ch_Z^X(V^{\bullet})=\ch_Z^U(V^{\bullet})\cdot[j].$$

\item[(iii)] Let $f:X'\ra X$ be a morphism, $Z'=f^{-1}(Z)$. Then
$$\ch_{Z'}^{X'}(f^*V^{\bullet})=f^*\ch_Z^X(V^{\bullet}).$$

\item[(iv)] If $0\ra V_1^{\bullet}\ra V^{\bullet}\ra V_2^{\bullet}\ra 0$
is an exact triple of $2$-periodic complexes exact off $Z$ then
$$\ch_Z^X(V^{\bullet})=\ch_Z^X(V_1^{\bullet})+\ch_Z^X(V_2^{\bullet}).$$

\item[(v)] If $W^{\bullet}$ is a $2$-periodic
complex of vector bundles on $X$
then
$$\ch_Z^X(W^{\bullet}\otimes V^{\bullet})=
\ch_Z^X(V^{\bullet})\cdot\ch(W^{\bullet}).$$

\item[(vi)]
Let $p:N\ra X$ be a vector bundle over $X$ and $\We^{-\bullet}p^*N^{\vee}$ be
the Koszul-Thom complex on $N$ (concentrated in degrees $[-\rk N,0]$)
which is exact off the zero section $i:X\hra N$.
Then one has
$$
\ch_Z^N(\We^{-\bullet}
p^*N^{\vee}\otimes p^*V^{\bullet})=
\td(N)^{-1}\cdot\ch_Z^X(V^{\bullet})\cdot [i].$$
\end{enumerate}
\end{prop}

To prove (iv) one has to use Lemma \ref{deform} to 
reduce to the case of a split exact triple. 
The proof of (vi) is similar to that of \cite[Prop.\ 3.4, Ch.\ 1]{BFM}.
\medskip

\begin{rems}

\

1. It seems natural to conjecture that for a pair
of $2$-periodic complexes of vector bundles
$V_1^{\bullet}$ and $V_2^{\bullet}$ such that $V_1^{\bullet}$
(resp. $V_2^{\bullet}$)
is exact outside of $Z_1\subset X$ (resp. $Z_2\subset X$) one has
$$\ch_{Z_1\cap Z_2}^X(V_1^{\bullet}\otimes V_2^{\bullet})=
\ch_{Z_1}^X(V_1^{\bullet})\cup\ch_{Z_2}^X(V_2^{\bullet}).$$
However, we did not check this property (cf.~\cite[Example 18.1.5]{F}).

2. For a $2$-periodic complex $V^{\bullet}$ 
of vector bundles on $X$ which
is exact off $Z$  we can consider the class
$[H^+(V^{\bullet})]-[H^-(V^{\bullet})]$ in $K_0(Z)$ (the Grothendieck
group of coherent sheaves on $Z$) defined by the difference
of even and odd cohomology of $V^{\bullet}$.
Applying the usual graph-construction
for bounded complexes we get the localized Chern character
$\ch_Z^X([H^+(V^{\bullet})]-[H^-(V^{\bullet})])$.
It would be interesting to check whether this class coincides
with $\ch_Z^X(V^{\bullet})$.
\end{rems}

\section{Orthogonal bundles with isotropic sections}

\subsection{Main construction}

We start with the following data. Let $E$ be an orthogonal bundle over a variety
$X$ of even rank $2n$ with an isotropic section $s:X\ra E$ and
$S=S^+\oplus S^-$ be a corresponding spinor bundle 
(an irreducible representation of the Clifford algebra of $E$)
with fixed decomposition into even and odd components.
Note that locally the choice of $S$ is unique up to an isomorphism and up to
interchanging $S^+$ and $S^-$. 

Let $Z(s)\subset X$ be the zero locus of $s$.
Consider the action of the section $s$ on $S$. This gives the
following infinite $2$-periodic complex:
$$S^{\bullet}:
\ldots\ra S^+\stackrel{s}{\ra} S^-\stackrel{s}{\ra} S^+\ra\ldots$$
with $S^+$ at the degree $0$ term. It is well-known that $S^{\bullet}$
is strictly exact outside $Z(s)$.
Therefore, we can consider the corresponding
localized Chern character:
$$\ch(E,s,S):=\ch_{Z(s)}^X(S^{\bullet})\in A^*(Z(s)\ra X)_{\Q}.$$
According to Proposition~\ref{properties}.(ii), it has the following
property: 
\begin{equation}\label{i*S}
i_*\ch(E,s,S)=\ch(S^+)-\ch(S^-),
\end{equation}
where $i:Z(s)\hra X$ is the natural embedding.

Denote by $\Pi S$ the spinor bundle obtained from $S$ by changing
the parity, i.e.
$$(\Pi S)^+:=S^- \text{\ and \ }(\Pi S)^-:=S^+.$$
 Then the corresponding localized
Chern characters are related by
\begin{equation}\label{chpar}
\ch(E,s,\Pi S)=-\ch(E,s,S).
\end{equation} 

\subsection{Compatibility}

The main property of the class $\ch(E,s,S)$
constructed above is its compatibility with reduction with respect to
isotropic subbundles. 

Let $M\subset E$ be an isotropic
subbundle with the following properties:

\begin{enumerate}
\item[(a)] $s\mod M^{\perp}$ is a regular section of $E/M^{\perp}$,

\item[(b)] for all $x\in X$ the inclusion $s(x)\in M$ implies $s(x)=0$.
\end{enumerate}

In this situation we can consider the locus $X'\sub X$ of zeroes of
the regular section $s\mod M^{\perp}$ with $\dim X' = \dim X-\rk M$. 
The restriction of $M^{\perp}/M$ to $X'$ is an orthogonal bundle $E'$ 
equipped with a section $s'$ induced by $s$. 
Furthermore, if $S$ is a spinor bundle for $E$ then
the coinvariants $S_M=S/MS$ form a spinor bundle
for $M^{\perp}/M$. Restricting it to $X'$ we get a spinor bundle $S'$
for $E'$ and we can apply the above construction to
the triple $(E',s',S')$.
Note that by property (b) we have $Z(s')=Z(s)\subset X'\subset X$.

\begin{thm}\label{classproperty} In the described situation we have
$$\ch(E,s,S)=\td(M^{\vee})^{-1}\cdot\ch(E',s',S')\cdot [i]$$
where $[i]\in A^{\rk M}(X'\ra X)$ is the canonical orientation of
the regular embedding $i:X'\rInto X$.
\end{thm} 

\Pf . Consider first the
situation when $i:X'\ra X$ is the embedding of the zero section 
into a vector bundle and $E$ decomposes into orthogonal direct sum
of $M^{\perp}/M$ and $M\oplus M^{\vee}$ in such a way that $s$ is
a sum of an isotropic section $s'$ of $M^{\perp}/M$ and of a regular section
$s_0$ of $M^{\vee}\subset M\oplus M^{\vee}$ with zero locus $X'$. 
In this case, we have isomorphisms of $\Z/2\Z$-graded complexes 
$$S\simeq 
\Hom(\We{^\bullet} 
M^{\vee}, S_M)\simeq
 \We{^{\bullet}}
 M\otimes S_M$$
where the differential on $S_M$ is given by the action of $s'$ while
the differential on 
$\We{^{\bullet}} M$
 is $i(s_0)$. 
It remains to apply Proposition \ref{properties} (vi).

To reduce the general case to this one we use the standard deformation of
the embedding $X'\hra X$ to the normal cone and then apply
Lemma \ref{deform}. Recall that the total
space $\XX$ of this deformation is an
open subset in the blow-up of $X\times\P^1$
along $X'\times\{\infty\}$ obtained by removing the proper transform
of the divisor $X\times\{\infty\}$ (which is isomorphic to
$\Bl_{X'}X$). Then the fiber $\XX_{\infty}$ is isomorphic to the
normal cone of the embedding $X'\hra X$.
We want to construct a compatible deformation of the orthogonal data
$(E,s,M,S)$ (near $\infty$).
Consider the direct sum of two orthogonal bundles
$$\wt{E}=E\oplus (M\oplus E/M^{\perp})$$
on $X$, where the orthogonal structure on $M\oplus E/M^{\perp}$ is induced by 
the natural pairing between $M$ and $E/M^{\perp}$.
Then we have a family of isotropic subbundles (isomorphic to $M$)
in $\wt{E}$ parametrized by $\P^1$:
$$\wt{M}_{(\la_0:\la_1)}\hra \wt{E}: m\mapsto (-\la_1m,\la_0m,0).$$
Let
$$E_{(\la_0:\la_1)}=\wt{M}_{(\la_0:\la_1)}^{\perp}/\wt{M}_{(\la_0,\la_1)}$$
be the corresponding reduced orthogonal bundle.
Note that the subbundle $\wt{M}_{(\la_0:\la_1)}^{\perp}\subset\wt{E}$
is defined by the equations
$$\la_0\ov{v}\equiv \la_1 v\mod M^{\perp}$$
where $(v,m,\ov{v})\in\wt{E}$.
{}From this one can easily see that
$$E_{\infty}\simeq M^{\perp}/M\oplus(M\oplus M^{\vee}),$$
while $E_t\simeq E$ for $t\neq\infty$.
Let $U\subset\P^1$ be a small open neighborhood of $\infty$ (actually one
can take $U$ to be $\P^1\setminus\{0\}$). Then the
family $(\wt{M}_t)$ can be considered as an isotropic subbundle
$\wt{\MM}$ of $p^*\wt{E}$, where $p:X\times U\ra X$ is the projection.
Hence, the family $(E_t)$ defines an orthogonal
bundle $\EE=\wt{\MM}^{\perp}/\wt{\MM}$ on $X\times U$.
For every $t\in U$ we have an inclusion
$$0\oplus M\oplus 0\subset \wt{M}_t^{\perp}$$
which induces the embedding $M\hra E_t$ as a subbundle. Let
$\MM\subset \EE$ be the corresponding isotropic subbundle.

Notice that we have an isotropic section
$$s_{(\la_0:\la_1)}=(\la_0 s,0,\la_1 s\mod M^{\perp})$$
of $\wt{\MM}^{\perp}(1)$ which we will consider as a section of
$\wt{\MM}^{\perp}$ using the trivialization of $\OO_{\P^1}(1)$ over $U$.
Let $\wt{s}$ be the corresponding isotropic section of $\EE$.
Let $\pi:\XX\ra X\times\P^1$ be the natural projection. Let us denote
$\XX_U=\pi^{-1}(X\times U)$. Then $\pi^*\EE$ is the orthogonal
bundle on $\XX_U$. Note that the section $\wt{s}$ vanishes over
$X'\times\{\infty\}$, hence, its pull-back to $\XX_U$ vanishes
over $E\cap\XX_U$ where $E\subset\Bl_{X'\times\{\infty\}}X\times\P^1$
is the exceptional divisor. By definition of $\XX$ we have the equality
$E\cap\XX_U=\pi^{-1}(X\times\{\infty\})$. Assuming that
$\la_1\neq 0$ on $U$, we can take $\la=\la_0/\la_1$ as 
the local parameter at
$\infty\in\P^1$. 
We have shown that the section $\pi^*\wt{s}$ 
is divisible by $\la$ and we define an isotropic section $s_U$ of $\EE$
by setting $$s_U=\frac{\pi^*\wt{s}}{\la}.$$
It is easy to see that the induced section of $E_{\infty}$ has the
required form $(s',0,s_0)$.

It remains to construct the family of spinor bundles over $U$.
We can start with
$$\wt{S}=\Hom(
\We{^\bullet}
M^{\vee}, S)\simeq
\We{^{\bullet}}
M\otimes S$$
as a spinor bundle for $\wt{E}$. Then $\SS=p^*\wt{S}_{\wt{M}}$
is a spinor bundle for $\EE$ and $(\pi^*\EE,s_U,\pi^*\MM,\pi^*\SS)$
is the required deformation.
\ed

\subsection{Homogeneous class} 
Assume that in the above situation we
have an isotropic subbundle $L\subset E$ of maximal rank $n$
(by analogy with the symplectic case
we call such a subbundle {\it Lagrangian}). Then the coinvariants
$S_L$ form a $\Z/2\Z$-graded 
bundle on $X$. We denote by $\ch(S_L)$ its $\Z/2\Z$-graded Chern
character. The following theorem 
shows that in this situation, the class $\ch(E,s,S)$ multiplied by some
standard factors is homogeneous of degree $n$.

\begin{thm}\label{purity} The class 
$$\ch(E,s,S)\cdot\td(L^{\vee})\cdot\ch(S_L)^{-1}\in A^*(Z(s)\ra X)_{\Q}$$
is concentrated in degree $n$.
\end{thm}

We will deduce this theorem from the following result.

\begin{thm}\label{purity2} 
Let $E$ be an orthogonal bundle of rank $2n$ on a scheme
$Y$, $S$ be a spinor bundle for $E$,
$L\subset E$ be a Lagrangian subbundle, $M\subset E$ be an isotropic
subbundle of rank $1$. Then the class
$$\ch(S_M)\cdot\td(M^{-1})^{-1}\cdot
\td(L^{\vee})\cdot\ch(S_L)^{-1}\in A^*(Y)_{\Q}$$
is concentrated in degree $n-1$.
\end{thm}

\Pf . First of all we can assume that 
$E$ is isomorphic as an orthogonal bundle to $L\oplus L^{\vee}$.
Indeed, the choices of a Lagrangian complement to $L$ form
an affine bundle over $Y$ so it suffices to prove the statement for the 
pull-back to this affine bundle. Next we claim that it suffices to prove
our statement in the case when $Y$ is smooth. 
Indeed, there exists an embedding $i:Y\ra Y'$ into a smooth
scheme $Y'$ and a bundle $L'$ on $Y'$ such that $L\simeq i^*L'$
(see Lemma 18.2 of \cite{F}). Let $\cal Q\ra Y'$ be the
bundle of projective quadrics associated with the
orthogonal bundle $L'\oplus (L')^{\vee}$ (i.e. $\cal Q$ is a closed
subscheme in $\P(L'\oplus (L')^{\vee})$ consisting of isotropic lines).
Then we have an embedding of $Y$ into $\cal Q$
such that our data is the pull-back of the similar data on $\cal Q$. 
Since $\cal Q$ is smooth this proves our claim.

Thus, we can assume that $Y$ is smooth and connected. Consider the following
composition:
$$f:M\ra E\ra L^{\vee}.$$
Let us consider two cases.

{\em Case 1: $f=0$.}
In this case $M$ is a subbundle in $L$. It is clear that
the class in question does not change if we tensor $S$ by a line bundle,
or switch the parity of $S$. Therefore, it does not depend on $S$ at all,
so we can choose any $S$ we like. Let us choose $S=\We^{\bullet} L$.
Then $S_L\simeq\OO$, while $S_M\simeq\We^{\bullet}(L/M)$, so we have
$$\ch(S_M)\cdot\td(M^{-1})^{-1}\cdot\td(L^{\vee})=
\ch(\We{^\bullet}(L/M))
\cdot\td((L/M)^{\vee})
=c_{n-1}((L/M)^{\vee}).$$ 

\medskip

{\em Case 2: $f\neq 0$.}
In this case $f$ gives a non-zero section 
$\OO\ra M^{-1}\otimes L^{\vee}$. Let $J$ be the image of the dual map
$M\otimes L\ra\OO$, and let $\pi:\wt{Y}\ra Y$ be a blow-up along
the subscheme defined by $J$. Then we have
an embedding $\pi^*:A^*(Y)\ra A^*(\wt{Y})$ so it suffices to prove
the statement for the pull-back of our data to $\wt{Y}$. But on $\wt{Y}$
the map $\pi^*f$ factors as follows:
$$\pi^*f:\pi^*M\ra \pi^*M(E)\hra L^{\vee}$$
where $E$ is the exceptional divisor, $N:=\pi^*M(E)$ embeds as
a subbundle into $L^{\vee}$. Let $N^{\perp}\subset L$ be the
orthogonal complement to $N$. We claim that $\pi^*M$ is contained
in $N^{\perp}\oplus N$. Indeed, it suffices to check this over the complement
to $E$, but there it follows from the condition that $M$ is isotropic.
Thus, we have a Lagrangian subbundle $L'=N^{\perp}\oplus N$ such
that $\pi^*M\subset L'$. Note that the statement of the theorem is true
with $\pi^*L$ replaced by $L'$ by Case 1 considered above.
It remains to prove that whenever we have
two Lagrangian subbundles $L_1$ and $L_2$ in $E$ such that
$L_1\cap L_2$ is a subbundle of rank $n-1$ in $E$ then the statements
of the theorem for $L_1$ and for $L_2$ are equivalent.
In fact, we claim that in this situation one has
$$\td(L_1^{\vee})\cdot\ch(S_{L_1})=-\td(L_2^{\vee})\cdot\ch(S_{L_2}).$$
Indeed, since $S_{L_1\cap L_2}$ is a spinor bundle for the rank-$2$
orthogonal bundle $L_1/(L_1\cap L_2)\oplus L_2/(L_1\cap L_2)$
everything reduces to the case $n=1$. In this case $L_2=L_1^{-1}$ and the
above statement follows from the identity
$$\td(L_1)^{-1}\cdot\td(L_1^{-1})\cdot\ch(L_1)=1.$$
\ed

\noindent
{\it Proof of Theorem \ref{purity}}.
Let us assume first that $s$ is a non-zero section of a line bundle 
$M$ which itself is embedded
as an isotropic subbundle into $E$. In this case by a deformation argument
(see the proof of Theorem~\ref{classproperty}) we can assume that
$E=M\oplus M^{-1}\oplus M^{\perp}/M$ as an orthogonal bundle. We have
$S\simeq \We^{\bullet}M\otimes S_M$, and therefore
$$\ch(E,s,S)=\ch(\OO\stackrel{s}{\ra}M)\cdot\ch(S_M).$$
Furthermore,
$$\ch(\OO\stackrel{s}{\ra}M)=-c_1(M,s)\cdot\td(M^{-1})^{-1}$$
where $c_1(M,s)\in A^1(Z(s)\ra X)$ is the localized first Chern class defined
by the section $s$ of $M$ (see \cite[Example 17.1.1]{F}).
It remains to apply Theorem~\ref{purity2} to conclude the proof.

To reduce to the case considered above we proceed as follows.
Let $p:C\cal Q\ra X$ be the affine quadrics fibration associated with $E$
(in other words, $C\cal Q$ consists of isotropic vectors in the total
space of $E$). Then we have a universal isotropic section of $s_0$ of
$p^*E$ and a morphism $f:X\ra C\cal Q$ such that $p\circ f=\id$ and
$s=f^*s_0$. By functoriality,
it suffices to prove our statement for the data $(p^*E,s_0,p^*S,p^*L)$.
Now let $\cal Q\ra X$ be the projective quadrics fibration associated with
$E$, $\wt{C\cal Q}\subset C\cal Q\times_X \cal Q$ be the incidence
correspondence. The natural projection $\pi:\wt{C\cal Q}\ra C\cal Q$ is
an isomorphism over the complement to the zero section $i:X\hra C\cal Q$,
while $\pi^{-1}(i(X))\simeq \cal Q$ is a locally trivial fibration over 
$i(X)$. It follows that every irreducible
variety in $C\cal Q$ is the birational image of some subvariety of 
$\wt{C\cal Q}$. According to \cite[Ex.\ 17.3.2]{F}, the pull-back map
$$\pi^*:A^*(Z(s_0)\ra C\cal Q)\ra A^*(Z(\pi^*s_0)\ra \wt{C\cal Q})$$
is injective, so it suffices to prove the statement for the pull-back
of our data by $\pi$. But
on $\wt{C\cal Q}$ our isotropic section $\pi^*s_0$ is contained in the 
isotropic subbundle of $\pi^*p^*E$ so we are done.
\ed

\section{Relative top Chern class}
\subsection{Definition on the level of complexes}\label{construction}

We refer to \cite{Del2} for the definition and main properties of the
symmetric powers of complexes of vector bundles. Below we will use the
fact that taking the $r$-th symmetric power (in characteristic zero)
is a well-defined functor on the homotopic category of
bounded complexes of vector bundles which preserves quasiisomorphisms.

Let 
$$C^{\bullet}=[C^0\stackrel{d}{\ra}C^1]$$ 
be a complex of vector bundles on a variety $X$
concentrated in degrees $0$ and $1$ equipped with a morphism of
complexes 
$$\tau:S^r(C^{\bullet})\ra\OO_X[-1].$$
 We assume that for every point $p\in X$ the induced pairing
\begin{equation}\label{paring}
S^{r-1}H^0(C^{\bullet}|_p)\otimes H^1(C^{\bullet}|_p)\ra\C
\end{equation}
has the following property: 

\begin{equation}\label{assumption}
\text{
{\em if $x^{r-1}$ is orthogonal to all of $H^1(C^{\bullet}|_p)$ for some
$x\in H^0(C^{\bullet}|_p)$ then $x=0$.}
}
\end{equation}

 Note that this condition automatically implies that $\rk C^0\le \rk C^1$. 
In this situation we are going to apply
the construction of the previous section to construct the 
class $c(C^{\bullet},\tau)\in A^*(X)_{\Q}$.

The complex $S^r(C^{\bullet})$ is as follows:
$$S^nC^0\rTo^{d_0} S^{r-1}C^0\otimes C^1\rTo^{d_1} 
S^{r-2}C^0\otimes \We{^2} 
C^1\ra\ldots$$
Hence, the morphism $\tau$ is given by the morphism
$$\tau_1:S^{r-1}C^0\otimes C^1\ra\OO_X$$
such that its composition with $d_0$ vanishes.
We can think about $\tau_1$ as of a section of a bundle on the total space of
$C^0$. Namely, if we denote by $p:C^0\ra X$ the natural projection then
$\tau_1$ gives rise to a section $\b\in H^0(C^0,p^*(C^1)^{\vee})$. 
On the other hand, 
the differential $d:C^0\ra C^1$ produces a section
$\a\in H^0(C^0,p^*C^1)$. 
It is easy to check 
that the pair of sections $(\a,\b)$ satisfies the following two conditions.

\begin{enumerate}
\item[(i)] the sections
$\a$ and $\b$ are orthogonal with respect to the natural
pairing 
$$H^0(C^0,p^*C^1)\otimes H^0(C^0,p^*(C^1)^{\vee})\ra 
H^0(C^0,\OO),$$

\item [(ii)] the intersection of the zero loci
$Z(\a)\cap Z(\b)$ coincides with the 
zero section  $X \subset C^0$.
\end{enumerate}

The bundle 
$E=p^*C^1\oplus p^*(C^1)^{\vee}$ on $C^0$ has a natural orthogonal structure,
and by  property (i) the section $s=(\a,\b)$ of $E$ is isotropic. 
As a spinor bundle we will take
$S=\We^{\bullet}p^*(C^1)^{\vee}$ with the natural grading (so that
$S^+=\oplus_n\We^{2n}p^*(C^1)^{\vee}$).
By property (ii) the zero locus $Z(s)$ coincides with
the zero section $X\subset C^0$.
Therefore, we have a well-defined class 
$$\ch(E,s,S)\in A(X\ra C^0)_{\Q}.$$

According to~\cite[Prop.\ 17.4.2]{F}, there is a natural isomorphism 
$$A(X\ra C^0)\ra A(X\ra X)=A^*(X):c\mapsto c\cdot [p]$$
where $[p]\in A^{-\rk C^0}(C^0\ra X)$ is the orientation class of $p$
and so we can define
\begin{equation}\label{cab}
c(C^{\bullet},\tau)=\td(C^1)\cdot\ch(E,s,S)\cdot[p]\in A^*(X)_{\Q}.
\end{equation}

Let $\chi(C^{\bullet})=\rk C^0-\rk C^1$ denote the Euler characteristic of
$C^{\bullet}$. 
{}From Theorem~\ref{purity} we immediately obtain that in fact
$$
c(C^{\bullet},\tau)\in A^{-\chi(C^{\bullet})}(X)_{\Q}.$$

We call $c(C^{\bullet},\tau)$ the {\it relative top Chern class of
$C^{\bullet}$ determined by}  $\tau$.
This definition is motivated by the following property.

\begin{prop}\label{difference} One has
$$c_{\rk C^0}(C^0)\cdot c(C^{\bullet},\tau)=c_{\rk C^1}(C^1).$$
\end{prop}

\Pf . Let $i:X\ra C^0$ be the zero section. Then by (\ref{i*S})
we have
$$i_*\td(C^1)\cdot\ch(E,s,S)=
\td(p^*C^1)\cdot\sum_n (-1)^n \ch(\We{^n} 
p^*(C^1)^{\vee})=c_{\rk C^1}(p^*C^1).$$
Now it is easy to see that for arbitrary class $c\in A^*(X\ra C^0)$ one has
$$i_*c=p^*(c_{\rk C^0}(C^0)\cdot c\cdot[p]).$$ 
Hence, we get
$$p^*(c_{\rk C^0}(C^0)\cdot c(C^{\bullet},\tau))=p^*(c_{\rk C^1}(C^1)).$$
Applying $i^*$ to both parts we get the result.
\ed

The first simple property of the relative top Chern class is the following.

\begin{prop}\label{homotopy} 
The class $c(C^{\bullet},\tau)$ depends only on the homotopy
class of $\tau$.
\end{prop}

\Pf . We will show that if $\tau_1$ is replaced by 
$\tau'_1=\tau_1+h\circ d_1$ for some map 
$$h:S^{r-2}C_0\otimes\We^2 C_1\ra\OO_X $$ 
then the corresponding
orthogonal data changes to an isomorphic one. Indeed,
$h$ gives rise to a section $\gamma\in H^0(C^0,\We^2p^*(C^1)^{\vee})$.
Thus, considering $\gamma$ as a skew-symmetric homomorphism
$p^* C^1\ra p^*(C^1)^{\vee}$ we can construct an automorphism of the
orthogonal bundle $p^*C^1\oplus p^*(C^1)^{\vee}$ acting as the identity
on the subbundle $p^*(C^1)^{\vee}$ and on the quotient by it.
It is easy to see that this automorphism sends the isotropic section
corresponding to $\tau_1$ to the isotropic section corresponding to $\tau'_1$
so we are done.
\ed

The main property of the relative top Chern class is its invariance
under quasiisomorphisms.

\begin{thm}\label{mainproperty} 
Let $f:\ov{C}^{\bullet}\ra C^{\bullet}$ be a quasiisomorphism
of complexes of bundles concentrated in degrees $0$ and $1$,
let $\tau:S^r(C^{\bullet})\ra\OO_M[-1]$ be a morphism satisfying the
assumption~(\ref{assumption}) above. Set $\ov{\tau}=\tau\circ S^r(f)$, then
$$c(\ov{C}^{\bullet},\ov{\tau})=c(C^{\bullet},\tau).$$
\end{thm}

Since by Proposition \ref{homotopy} we can work in the homotopic category
of complexes, we can use the following simple lemma.

\begin{lem} In the homotopic category
every quasiisomorphism of complexes of bundles concentrated in
degrees $0$ and $1$ is a composition of a quasiisomorphism 
which is an embedding of subbundles on each term followed by a
quasiisomorphism of the type $C^{\bullet}\oplus (A\stackrel{\id}{\ra} A)\ra
C^{\bullet}$.
\end{lem} 

\Pf . Let $f:\ov{C}^{\bullet}\ra C^{\bullet}$ be such a quasiisomorphism. 
We have a decomposition of $f$ into a composition of two quasiisomorphisms:
$f=f_1\circ f_2$ where
$$f_1:(\ov{C}^1\rTo^{\id}\ov{C}^1)\oplus C^{\bullet}\ra C^{\bullet}$$
is the natural projection,
$$f_2:\ov{C}^{\bullet}\ra (\ov{C}^1\rTo^{\id}\ov{C}^1)\oplus 
C^{\bullet}$$
is the composition of $f$ with the natural embedding of $C^{\bullet}$
into the target. Now we can replace $f_2$ by a homotopic map
$f'_2$ using the natural embedding 
$\ov{C}_1\ra \ov{C}_1\oplus C_0$ as a homotopy.
Now the fact that $f'_2$ is an embedding of subbundles on each term follows
from the exactness of the following sequence:
$$0\ra \ov{C}^0\ra \ov{C}_1\oplus C_0\ra C_1\ra 0.$$
\ed

\noindent
{\it Proof of Theorem \ref{mainproperty}}.
According to the above lemma it suffices to consider separately
two classes of quasiisomorphisms. 
Consider first the case when
$f_0:\ov{C}^0\hookrightarrow C^0$ and $f_1:\ov{C}^1\hookrightarrow C^1$ 
are embeddings of subbundles.
Let $\ov{p}:\ov{C}^0\ra X$ be the natural projection,
$\ov{s}=(\ov{\a},\ov{\b})$ be the isotropic section of 
$\ov{p}^*\ov{C}^1\oplus\ov{p}^*(\ov{C}^1)^{\vee}$ constructed from
the data $(\ov{C}^{\bullet},\ov{\tau})$.
It is easy to check the following two properties:

\begin{enumerate}
\item[(1)] $\ov{\b}$ is the image of 
$\b|_{\ov{C}^0}\in H^0(\ov{C}^0,\ov{p}^*(C^1)^{\vee})$ under 
the natural map 
$$H^0(\ov{C}^0,\ov{p}^*(C^1)^{\vee})\ra 
H^0(\ov{C}^0,\ov{p}^*(\ov{C}^1)^{\vee});$$

\item[(2)] the square
\begin{equation}
\begin{diagram}
\ov{C}^0 & \rTo & C^0\\
\dTo   & & \dTo^{\a}\\
p^*\ov{C}^1 &\rTo^i & p^*C^1
\end{diagram}
\end{equation}
is Cartesian (and in particular, commutative), where the left vertical
arrow is the composition of $\ov{\a}$ 
with the natural embedding $\ov{p}^*\ov{C}^1\hra p^*\ov{C}^1$.
\end{enumerate}

The condition (2) says, in particular, that $\ov{C}^0$ is the zero locus
of the regular section $\a\mod p^*\ov{C}^1$ of $p^*C^1/p^*\ov{C}^1$.

Let $(\ov{C}^1)^{\perp}\subset (C^1)^{\vee}$ be the subbundle
orthogonal to $\ov{C}^1\subset C^1$.
Then $M:=p^*(\ov{C}^1)^{\perp}$ is an isotropic subbundle in 
$E=p^*C^1\oplus p^*(C^1)^{\vee}$. 
We claim that $M$ satisfies the conditions (a) and (b) preceding
Theorem~\ref{classproperty} and that the reduced data correspond
to $(\ov{\a},\ov{\b})$. Indeed, the orthogonal complement of $M$ is $E$ is
$$M^{\perp}=p^*\ov{C}^1\oplus p^*(C^1)^{\vee}.$$
Thus, $s\mod M^{\perp}$ is just the section $\a \mod p^*\ov{C}^1$. 
This is a regular section of $E/M^{\perp}$ and its zero locus is
$\ov{C}^0\subset C^0$ (in particular, the condition (a) holds). 
Note that we have a natural identification
$$M^{\perp}/M\simeq p^*\ov{C}^1\oplus p^*(\ov{C}^1)^{\vee}.$$
It is easy to see that under this identification
the section of $E'=M^{\perp}/M|_{\ov{C}^0}$ induced
by $s=(\a,\b)$ coincides with $s'=(\ov{\a},\ov{\b})$. In particular,
the zero locus of $\ov{s}$ coincides with the zero section of $\ov{C}^0$,
hence the condition (b) holds.
Finally, the coinvariants of $M=p^*(\ov{C}^1)^{\perp}$ in the spinor bundle
$S=\We^{\bullet}p^*(C^1)^{\vee}$ for $E=p^*(C^1)\oplus p^*(C^1)^{\vee}$ 
give the
spinor bundle $\We^{\bullet}p^*(\ov{C}^1)^{\vee}$ for $M^{\perp}/M$ which
restricts to the spinor bundle
$S'=\We^{\bullet}\ov{p}^*(\ov{C}^1)^{\vee}$ for $E'$. 
Now Theorem \ref{classproperty} gives the equality
$$\ch(E,s,S)=\td(C^1/\ov{C}^1)^{-1}\cdot\ch(E',s',S')\cdot [i]$$
where $i:\ov{C}^0\ra C^0$ is the natural embedding.
Hence,
\begin{eqnarray*}
c(C^{\bullet},\tau)&=&\td(C^1)\cdot\ch(E,s,S)\cdot[p]=
\td(C^1)\cdot\td(C^1/\ov{C}^1)^{-1}\cdot\ch(E',s',S')\cdot [i]\cdot [p]\\
&=&
\td(\ov{C}^1)\cdot\ch(E',s',S')\cdot [\ov{p}]=c(\ov{C}^{\bullet},\ov{\tau}).
\end{eqnarray*}

The case of the quasiisomorphism $C^{\bullet}\oplus (A\rTo^{\id} A)
\ra C^{\bullet}$ is much easier and is left to the reader.
\ed

Proposition~\ref{difference} admits the following generalization.

\begin{prop}\label{subbundle}
Let $\ov{C}^{\bullet}\subset C^{\bullet}$ be a subcomplex of $C^{\bullet}$
such that $\ov{C}^0\subset C^0$ is a subbundle and $\ov{C}^1=C^1$.
Then we have
\begin{equation}\label{chern}
c(\ov{C}^{\bullet},\ov{\tau})=
c_{top}(C^0/\ov{C}^0)\cdot c(C^{\bullet},\tau)
\end{equation}
\end{prop}

\Pf . Let $i:\ov{C}^0\ra C^0$ be the embedding of the total spaces
of vector bundles. Then by Proposition~\ref{properties}.(iii)
we have
$$i^*\ch(E,s,S)=\ch(i^*E,i^*s,i^*S),$$
where $(E,s,S)$ are the orthogonal data on $C^0$ constructed from
the data $(C^{\bullet},\tau)$. Now our statement follows from the fact
that for arbitrary class $c\in A^*(X\ra C^0)$ one has
$$i^*(c)\cdot[p\circ i]= c_{top}(C^0/\ov{C}^0)\cdot c\cdot[p].$$
This can be easily proved by reducing to the case when the
embedding $\ov{C}^0\hra C^0$ splits.
\ed

Proposition \ref{difference} corresponds to the case $\ov{C}^0=0$.

\subsection{Passage to the derived category}
\label{derived}

Now assume that we have a complex $C^0\ra C^1$ as before
and a map $\tau:S^r(C^{\bullet})\ra\OO_X[-1]$ {\it in the derived category}
$D^b(X)$ such that the condition~(\ref{assumption}) on the induced
cohomological pairing~(\ref{paring})
is satisfied. Assuming that $X$ is quasi-projective we will show that one can
still define the relative top Chern class $c(C^{\bullet},\tau)$ in this
situation. Let $\OO_X(1)$ be an ample line bundle on $X$. 
Let us denote by $K^b(X)$ the homotopic category of bounded
complexes of coherent sheaves on $X$.

The proof of the following lemma is straightforward.

\begin{lem}\label{straight} Every surjection $f^1:\ov{C}^1\rOnto C^1$ 
extends uniquely to a quasiisomorphism
$$f:[\ov{C}^0\ra\ov{C}^1]\ra [C^0\ra C^1].$$
\end{lem}

\begin{prop}\label{resolve} 
In the above situation for a sufficiently large $m$ 
and a surjection of the form 
$$\OO_X(-m)^{\oplus N}\rOnto C^1,$$
let $\ov{C}^{\bullet}\ra C^{\bullet}$ be the corresponding 
quasiisomorphism, where $\ov{C}^1=\OO_X(-m)^{\oplus N}$. Then
$$\Hom_{K^b(X)}(S^r(\ov{C}^{\bullet}),\OO_X[-1])\simeq
\Hom_{D^b(X)}(S^r(\ov{C}^{\bullet}),\OO_X[-1]).$$
\end{prop}

\Pf . 
Denote $S^r(\ov{C}^{\bullet})$ by $E^{\bullet}$.
{}From the spectral sequence that computes
$\Hom_{D^b(X)}(E^{\bullet},\OO_X[-1])$, it is clear
that it is sufficient to prove that $\Ext^i(E^j,\OO_X)=0$
for $i>0$, $j>0$. Since $E^j=S^{r-j}(\ov{C}^0)\otimes\We^j(\ov{C}^1)$,
it is enough to prove that
$H^i(S^j(\ov{C}^0)^{\vee}(m))=0$ for $i>0$, $j<r$.
We know that for all sufficiently
large $m$ one has
\begin{equation}\label{vanish}
H^{>0}(S^j(C^0)^{\vee}\otimes S^{j_1}(C^1)^{\vee}\otimes\ldots\otimes
S^{j_k}(C^1)^{\vee}(m))=0
\end{equation}
as long as $j+j_1+\ldots+j_k<r$.
Let $K$ be the kernel of the surjection $\ov{C}^i\ra C^i$ ($i=1,2$).
Then from the exact sequence
$$0\ra (C^0)^{\vee}\ra (\ov{C}^0)^{\vee}\ra K^{\vee}\ra 0$$
we see that vanishing of $H^{>0}(S^j(\ov{C}^0)^{\vee}(m))$
follows from vanishing of 
$$H^{>0}(S^j(C^0)^{\vee}\otimes S^{j_1}K^{\vee}(m))$$ 
for $j+j_1<r$. On the other hand, we have
the exact sequence
$$0\ra (C^1)^{\vee}\ra (\ov{C}^1)^{\vee}\ra K^{\vee}\ra 0.$$
Considering the $j_1$-th symmetric power of this
sequence and using (\ref{vanish}) we see that it is sufficient to prove 
vanishing of
$H^{>0}(S^j(C^0)^{\vee}\otimes S^{j_1}(C^1)^{\vee}\otimes
S^{j_2}K^{\vee}(m))$ for $j+j_1+j_2<r$, $j_1>0$. Now we can continue arguing
in the same way using (\ref{vanish}).
\ed

\begin{lem}\label{ext1} 
Let $C^0\ra C^1$ and $D^0\ra D^1$ be two complexes of sheaves
on $X$ such that $\Ext^1(C^1,D^0)=0$. Then the natural map
$$\Hom_{K^b(X)}(C^{\bullet},D^{\bullet})\ra
\Hom_{D^b(X)}(C^{\bullet},D^{\bullet})$$
is an isomorphism.
\end{lem}

\medskip

The class $c(C^{\bullet},\tau)$ can now be defined as follows. 

Using Proposition~\ref{resolve} we can choose a quasiisomorphism
$\ov{C}^{\bullet}\ra C^{\bullet}$ such that $\tau$ is realized by
a map of complexes 
$$\ov{\tau}:S^r(\ov{C}^{\bullet})\ra\OO_X[-1]$$
and apply the construction of 
Section~\ref{construction} to $\ov{C}^{\bullet}$. 
We set
$$ c(C^{\bullet},\tau):=c(\ov{C}^{\bullet},\ov{\tau}).$$

To check that this is well-defined,
assume that we have another quasiisomorphism
$\wt{C}^{\bullet}\ra C^{\bullet}$ such that $\tau$ is realized
by a map of complexes 
$$\wt{\tau}:S^r(\wt{C}^{\bullet})\ra\OO_X[-1].$$
Then we can find a quasiisomorphism
$D^{\bullet}\ra\ov{C}^{\bullet}$ with $D^1=\OO_X(-m)^{\oplus N}$
where $m$ is sufficiently large. Then by Lemma~\ref{ext1}
the following composition of morphisms in the derived category
$$D^{\bullet}\ra\ov{C}^{\bullet}\ra C^{\bullet}\ra\wt{C}^{\bullet}$$
is realized by a map of complexes.

Now the equality $c(\ov{C}^{\bullet},\ov{\tau})=c(\wt{C}^{\bullet},\wt{\tau})$
follows from Theorem~\ref{mainproperty} and Propositions~\ref{homotopy}
and~\ref{resolve}.

\section{Witten's top Chern class}

\subsection{Construction of the spin virtual class}
\label{mainconstr}

Witten's ``top Chern class'' is defined in the following situation.

Let $\pi:\CC\ra S$ be a flat family of prestable 
curves over a quasi-projective base $S$,
and $\TT$ be a family of rank-one torsion-free sheaves on $\CC$ with
a non-zero homomorphism
$$b:\TT^r\ra\om_{\CC/S},$$
 where $\om_{\CC/S}$ is the relative
dualizing sheaf of $\pi$. 
(This data is part of an $r$-spin  structure on the family $\CC\ra S$;
see \cite{J2,J,JKV} for details.)
  
We can represent the sheaf $R\pi_*(\TT)$ by a complex 
$$C^{\bullet}=[C^0\ra C^1],$$
where $C^i$ are vector bundles on $S$. 
For example, we can choose an ample relative Cartier divisor $D\subset
\CC$ and consider the resolution of $\TT$ by the complex
$\TT(ND)\ra \TT(ND)/\TT$ where $N$ is large enough. Since the terms of this
complex are acyclic for $R\pi_*$, we can take $C^{\bullet}$ to be
$$\pi_*(\TT(ND))\ra \pi_*(\TT(ND)/\TT).$$

The morphism $b$ induces a morphism 
$$S^r(C^{\bullet})\ra R\pi_*(\om_{\CC/S})$$
in the derived category. 
Composing it with the canonical trace morphism 
$$R\pi_*(\om_{\CC/S})\ra\OO_S[-1],$$ 
we obtain
a morphism (still in the derived category) 
$$\tau:S^r(C^{\bullet})\ra\OO_S[-1].$$
We claim that for every $s\in S$, the induced pairing
$$S^{r-1}H^0(\CC_s,\TT_s)\otimes H^1(\CC_s,\TT_s)\ra\C$$
satisfies the condition~(\ref{assumption}).
Indeed, this pairing corresponds to the following composition
$$S^{r-1}H^0(\CC_s,\TT_s)\ra H^0(\CC_s,\TT^{r-1}_s)\ra
\Hom(\TT,\om_{\CC_s})\simeq H^1(\CC_s,\TT_s)^*$$
so the required condition follows from the injectivity of the second
arrow. 

Now we define the class
$$c(\TT,b):=c(\CC^{\bullet},\tau)\in A^*(S)_{\Q}$$
by applying the construction of
Section~\ref{derived}. The results of that section imply that
it does not depend on a choice of a complex $C^{\bullet}$
representing $R\pi_*(\TT)$.

\subsection{Spin virtual class}

In~\cite{JKV} axioms have been formulated for a so-called
spin virtual class $\cv$ (or Witten's top Chern class)
on the moduli space of $r$-spin curves so that it can be
used to produce a cohomological quantum fields theory.
Since every family of $n$-punctured $r$-spin curves $\CC\to S$ of
{\em non-negative} type $(m_1,\ldots,m_n)$ (i.e.\ $m_j \ge 0$ for
$j=1,\ldots,n$) contains in its structural data a flat family of torsion-free
rank-one sheaves $\TT$ on $\CC$ with a map 
$$b:\TT^r \to \om_{C/\SS},$$
the above construction produces a Chow cohomology class
$c(\TT,b)$ on $S$ of degree $-\chi(\TT)$. We propose this class
as a candidate for the spin virtual class. In the remainder of the
paper, we will consider examples and establish some properties
of this class which justify our suggestion.

\subsection{Examples}

\

 1. Assume that the spin structure is {\em convex},
i.e.\ that $R^0\pi_*(\TT)=0$ and therefore $R^1\pi_*(\TT)$ is
a bundle. 
Then from Proposition~\ref{difference} it follows
$$c(\TT,b)=c_{top}(R^1\pi_*(\TT)),$$
and so in this situation our class is a genuine top Chern class
of a bundle.
\medskip

 2. Assume that $r=2$, the map $\pi$ is smooth, $S$ is connected, 
$\TT$ is locally free and the map
$b:\TT^2\ra\om_{\CC/S}$ is an isomorphism. In other words,
$\CC\to S$ is a family of smooth curves with ordinary spin structures
(theta characteristics).

Then $\chi=0$, so the class $c(\TT,b)$ has degree $0$ i.e.\ it
is just a number.
We claim that this number is equal to $(-1)^{h^0(\CC_s,\TT_s)}$
for any $s\in S$. Indeed, we can assume that $S$ is a point, $\CC=C$
is a smooth curve. Then we have to compute the class corresponding to 
the complex $H^0(C,\TT)\ra H^1(C,\TT)$ with zero differential and
to the map 
$$H^0(C,\TT)\ra H^1(C,\TT)^*$$ 
which is an isomorphism
obtained from Serre duality. Now our statement follows from
Proposition~\ref{properties}(vi). The sign appears because the complex
arising from the action of the isotropic section on the spinor bundle  
coincides with the Koszul complex shifted by $h^0(\CC_s,\TT_s)$.
\medskip

 3. Let $\CC \to S$ be a family of elliptic curves
with $\TT$ trivial on every fiber $\CC_s$ such that the map
$b:\TT^r\ra\om_{\CC/S}$ is an isomorphism. Then we claim that
$c(\TT,b)=-r+1$. 

Indeed, again the class $c(\TT,b)$ has degree $0$
and we can assume that $S$ is a point.
Now we have to compute the class corresponding to the complex
$C^0\ra C^1$ with zero differential where $C^0=C^1=\C$ and the map
$S^{r-1}C^0\ra (C^1)^*$ is an isomorphism. Thus, we are reduced to computing
the localized Chern class of the complex 
$$\OO_{\A^1}\rTo^{z^{r-1}}\OO_{\A^1}$$
concentrated in degrees $[0,1]$, where $z$ is the coordinate on $\A^1$. 
It is easy to see that it is equal to $-(r-1)$.

\subsection{Factorization properties}

Let $\Gamma$ be a stable decorated genus $g$ graph that corresponds to
a stratum $\mbar_\Gamma^{1/r}$ in the moduli stack
$\mbar_{g,n}^{1/r}$ of $r$-spin genus $g$ curves. (We adopt the
notation from~\cite{JKV}.) 
Denote by $\wt{\Gamma}$ the decorated graph obtained by cutting all
edges of $\Gamma$ and by $\Gamma_0$ and $\wt{\Gamma}_0$ the corresponding
underlying graphs (without decorations). Let $\mbar_{\Gamma_0}$
and $\mbar_{\wt{\Gamma}_0}$ be the corresponding moduli spaces of
curves. For each of these moduli spaces $\mbar$, let $\CC \to \mbar$
be the corresponding universal curve. Finally, consider the fibered
product
$$
X=\mbar_{\wt{\Gamma}_0} \times _{\mbar_{{\Gamma}_0}} \mbar_{\Gamma}^{1/r}
$$
and the pull-backs $\CC_X$ and $\wt{\CC}_X$ of the universal curves
$\CC_{\Gamma}$ and $\CC_{\wt{\Gamma}}$ to $X$.

Consider the following diagram
$$
 \begin{diagram} 
&&\CC_{\wt{\Gamma}}&\lTo^{\wt{p}} & \wt{\CC}_X& & \\
&&   &  & \dTo_\nu & & \\
&&  \dTo_{\wt{\rho}} & & \CC_X & & \\
&&   &  &  &\rdTo_{\hat{q}} & \\
&&   &  &\dTo_{\pi_X}  & &\CC_\Gamma \\
\mbar_{\wt{\Gamma}_0}&\lTo &\mbar_{\wt{\Gamma}}^{1/r} &\lTo^p &X &&\dTo^\pi\\ 
&\rdTo  &  &  & &\rdTo_q & \\
& &\mbar_{\Gamma_0}&&\lTo& &\mbar_{\Gamma}^{1/r} 
 \end{diagram}
$$
where $\nu$ is the normalization morphism and the morphism $p$
is given by the canonical construction that induces an $r$-spin structure
on the normalization of a spin curve (see~\cite{J2}).
Let $\TT$ and $\wt{\TT}$
be the universal $r$-th root sheaves on $\CC_{\Gamma}$ and
$\CC_{\wt{\Gamma}}$, then by~\cite{J2} we have
$$\nu_*\wt{p}^*\wt{\TT} = \hat{q}^* \TT. $$
Since all the quadrilaterals in the diagram above are Cartesian,
we can compare the virtual top Chern classes corresponding
to $\TT$ and $\wt{\TT}$.

 Indeed, we have
$$q^*R\pi_*\TT=R\pi_{X*}\hat{q}^*\TT=R\pi_{X*}\nu_*\wt{p}^*\wt{\TT}=
R(\pi_X\nu)_*\wt{p}^*\wt{\TT}=p^*R\wt{\rho}_*\wt{\TT}.$$

Since both $p$ and $q$ are finite flat
morphisms and the virtual top Chern class commutes with \'etale 
base change, the above equation gives
$$q^*c(\TT,b) = p^* c(\wt{\TT},\wt{b})$$
for all stable decorated graphs $\Gamma$ except those that have a
non-separating  
edge marked by the pair $(r-1,r-1)$ (i.e.\ the spin structure $\TT$ at the
corresponding node is Ramond). For such $\Gamma$ the induced spin structure
$\wt{\TT}$ will not be of a non-negative type and the map $\wt{b}$ 
(and therefore the class $c(\wt{\TT},\wt{b})$) will not be defined.

This implies that the class $c(\TT,b)$ satisfies
Axiom~1 of the spin virtial class of~\cite{JKV}
and also Axiom~3 when $\Gamma$ has no separtating edges marked
with $(r-1,r-1)$. 
\medskip

Since Axioms~2 (Convexity) and~5 (Forgetting tails)
are obviously satisfied as well, it only remains to
verify Axiom~4 and also Axiom~3 in the case  when $\Gamma$ has
a non-separating edge corresponding to a Ramond singularity.
Axiom~4 (Vanishing) requires that $c(\TT,b)=0$ if $\Gamma$ has a tail marked
by $r-1$, and modulo this axiom the remaining case of Axiom~3 can  also be
stated as a vanishing condition that $p_*q^*c(\TT,b) = 0$.

In this paper we will only check the vanishing for $r=2$.

\subsection{Vanishing axiom in the case of square roots}

Assume that we have a family of curves $\pi:\CC\ra S$, a collection of 
sections $p_1,\ldots,p_{2k}:S\ra\CC$, and a family of torsion free rank-$1$
sheaves $\TT$ which is locally free near $p_i$'s
together with a morphism
$$b:\TT^2\ra\om_{\CC/S}(-p_1-\ldots-p_{2k})$$
 which induces an isomorphism
\begin{equation}\label{isomor}
\TT \isomarrow R\underline{\Hom}(\TT,\om_{\CC/S})(-p_1-\ldots-p_{2k})
\end{equation}
(see~\cite[Note 2.1.5]{J2}).

\medskip

Let $C^0\ra C^1$ be a complex of vector bundles representing
$R\pi_*(\TT)$ such that the map $S^2 R\pi_*(\TT)\ra\OO_S[-1]$
is realized on the level of complexes. Then we have a natural
map $C^0\ra (C^1)^{\vee}$ such that the map
$$C^0\ra C^1\oplus (C^1)^{\vee}$$ is an embedding of an isotropic
subbundle.
Note that the natural morphism of complexes $C^{\bullet}\ra 
(C^{\bullet})^{\vee}[-1]$ represents the natural map
$$R\pi_*(\TT)\ra R\pi_*R\underline{\Hom}(\TT,\om_{\CC/S})$$
in $D^b(S)$ induced by $b$.
Its cone shifted by $[-1]$ is the complex
\begin{equation}\label{cone}
C^0\ra C^1\oplus (C^1)^{\vee}\ra (C^0)^{\vee}.
\end{equation}
On the other hand, from the isomorphism (\ref{isomor}) it follows
that the cone of the above map in $D^b(S)$ is isomorphic to 
$\TT(p_1+\ldots+p_{2k})|_{p_1+\ldots+p_{2k}}$. Now we notice that
the complex (\ref{cone}) is quasiisomorphic to the complex
$$C^0\ra (C^0)^{\perp},$$
where $(C^0)^{\perp}$ is the orthogonal complement to $C^0\subset C^1\oplus
(C^1)^{\vee}$.
 Hence, we have the following isomorphism of vector bundles on
$S$: 
$$(C^0)^{\perp}/C^0\simeq \bigoplus_i \TT(p_i)|_{p_i}.$$
Furthermore, this isomorphism is compatible with the orthogonal
structures on both bundles where on the second bundle
we have the  
direct sum of the orthogonal structures on the line bundles
$\TT(p_i)|_{p_i}$ given by the isomorphism
$$(\TT(p_i)|_{p_i})^{\otimes 2}\ra\om(p_i)|_{p_i}\simeq\OO_S$$
(induced by $b$). Consider a finite \'etale covering $S'\ra S$
which makes all line bundles $\TT(p_i)|_{p_i}$ trivial.
Then on this covering the orthogonal bundle $(C^0)^{\perp}/C^0$ can
be trivialized (as an orthogonal bundle). In particular, we can choose
a trivial Lagrangian subbundle $\ov{M}\subset (C^0)^{\perp}/C^0$.
Let $M\subset (C^0)^{\perp}$ be the corresponding Lagrangian
subbundle in $C^1\oplus (C^1)^{\vee}$. Then we have an embedding of complexes
$$[C^0\ra C^1]\ra [M\ra C^1]$$ such that the morphism
$S^2(C^0\ra C^1)\ra\OO_S[-1]$ is induced by the natural morphism
$S^2(M\ra C^1)\ra\OO_S[-1]$ (which uses the map $M\ra (C^1)^{\vee}$).
Hence, according to Proposition \ref{subbundle}, the relative top Chern class
of $C^0\ra C^1$ is equal to the product of the relative top Chern class
of $M\ra C^1$ with $c_{top}(M/C^0)$. But $M/C^0\simeq\ov{M}$ is a trivial
bundle, hence $c_{top}(M/C^0)=0$. 

\subsection{Descent axiom}

\

Here we show that our virtual top Chern class satisfies the 
Descent Axiom formulated in~\cite{JKV2}. It is shown in~\cite{JKV2}
that this axiom implies the Vanishing axiom provided the virtual
class is defined for families of spin structures of type
$(m_1,\ldots,m_n)$, where one of $m_j$ can be equal to $-1$ (and
the rest are non-negative). Our construction, however, works only
for spin structures with $m_j \ge 0$ for all $j$ and therefore 
the question whether the Vanishing Axiom holds for the class $c(\TT,b)$ 
remains open.

Let $\si:S\ra\CC$ be a section such that 
$\pi$ is smooth along $\si(S)$ and $\TT$ is locally free along $\si(S)$. 
We can equip the sheaf
$\ov{\TT}=\TT(-\si)$ with the map
$$\ov{b}: \ov{\TT}^r\ra\om_{\CC/S}$$
induced by $b$.

\begin{prop} In the above situation one has 
$$c(\ov{\TT},\ov{b})=c(\TT,b)\cdot c_1(\si^*\TT).$$
\end{prop}

\Pf . {}From the exact sequence
$$0\ra \ov{\TT}\ra\TT\ra\TT|_{\si(S)}\ra 0$$
we obtain the exact triangle 
$$R\pi_*(\ov{\TT})\ra R\pi_*(\TT)\ra L$$ 
in the derived category $D^b(S)$, where $L:=\si^*\TT$.
Let us choose a complex of vector bundles $C^{\bullet}=[C^0\ra C^1]$
representing $R\pi_*(\TT)$ in such a way that $\Ext^1(C^1,L)=0$
(this is possible by Lemma~\ref{straight}). 
Then by Lemma \ref{ext1} the morphism
$(C^{\bullet})\ra L$ in $D^b(S)$ is represented by a morphism
of complexes, i.e. by a map $C^0\ra L$. Let us replace
$C^{\bullet}$ by the quasiisomorphic complex
$C^{\bullet}\oplus (L\rTo^{\id}L)$
and replace the induced map
$C^{\bullet}\oplus (L\rTo^{\id}L)\ra L$
by the homotopic one using the projection $C^1\oplus L\ra L$ as a homotopy.
By this procedure we represent the map $R\pi_*(\TT)\ra L$
by a morphism of complexes $f:C^{\bullet}\ra L$ such that $C^0\ra L$ is
surjective. Let $\ov{C}^{\bullet}$ be a kernel of $f$. Then
$\ov{C}^{\bullet}$ is a complex representing
$R\pi_*(\ov{\TT})$. It remains to apply Proposition \ref{subbundle}.
\ed

\subsection{Geometric interpretation}

In this section we will consider some geometric consequences of
the vanishing of the virtual top Chern class in the case of square 
roots of $\om(-p_1-p_2)$. Let us fix a smooth projective curve 
$C$ of genus $g$. 
Let $\pi:S\ra S^2C-\Delta$ be an \'etale covering of degree $2^{2g}$ 
corresponding to
a choice of a square root of $\om_C(-p_1-p_2)$ for $(p_1,p_2)\in S^2C-\Delta$.
Then we have a line bundle $\TT$ on $S\times C$ such that
$\TT^2\simeq p_C^*\om_C(-\DD)$ where $\DD\subset S\times C$ is the universal
divisor, $p_C:S\times C\ra C$ is the projection. 
Let $Z\subset S$ be the locus of $s\in S$ such that 
$H^0(C,\TT_s)\neq 0$. Then over $S - Z$ the virtual top Chern class is just
the first 
Chern class of the line bundle $R^1p_{S*}(\TT)$ where $p_S:S\times C\ra S$
is the projection. It is easy to show directly that the square of this
line bundle is trivial (see the proof of the theorem below). We are 
going to use Deligne's theory of determinant bundles in order to 
interpret this fact geometrically. 

\begin{thm} For a smooth projective curve $C$ of genus $g$ let
$X\subset S^2C-\Delta$ be the closed subset consisting of pairs
$(p_1,p_2)$ such that there exists an effective divisor $D$ with
$\OO_C(2D+p_1+p_2)\simeq\om_C$. Let
$U=S^2C-\Delta-X$ be the open complement to $X$. Then the line bundle
$\om_U^{\otimes 2^{2g}}$ on $U$ is trivial.
\end{thm}

\Pf . First of all let us explain why on the locus where 
$H^0(\TT)=0$, one has a trivialization of $H^1(\TT)^{\otimes 2}$.
Consider the exact sequence
$$0\ra H^0(\TT(p_1))\ra\TT(p_1)|_{p_1}\ra
H^1(\TT)\ra\ldots$$
We claim that the boundary homomorphism $\delta:\TT(p_1)|_{p_1}\ra H^1(\TT)$
is an isomorphism. Indeed, assume it is not. Then it is equal to zero,
hence the restriction map 
$H^0(\TT(p_1))\ra\TT(p_1)|_{p_1}$ is an isomorphism.
On the other hand, by restricting via $\delta$ 
the canonical pairing between
$H^1(\TT)$ and $H^0(\om\otimes\TT^{-1})$ to a pairing
between $\TT(p_1)|_{p_1}$ and $H^0(\om\otimes\TT^{-1})$ we get
the natural pairing
$$\TT(p_1)|_{p_1}\otimes H^0(\om\otimes\TT^{-1})\ra
\TT(p_1)|_{p_1}\otimes (\om\otimes\TT^{-1})|_{p_1}\ra\om(p_1)|_{p_1}\simeq
\OO_{p_1}.$$
The fact that this pairing is zero implies that the restriction
map $$H^0(\om\otimes\TT^{-1})\ra (\om\otimes\TT^{-1})|_{p_1}$$
 is zero.
But we have $\om\otimes\TT^{-1}\simeq\TT(p_1+p_2)$, and we know
that on the subspace $H^0(\TT(p_1))$ the above restriction map is not zero.
The obtained contradiction proves that we have an isomorphism
$$\TT(p_1)|_{p_1}\wt{\ra} H^1(\TT).$$
On the other hand, $(\TT(p_1)|_{p_1})^{\otimes 2}$ is trivial, hence,
$H^1(\TT)^{\otimes 2}$ is trivial.
Since $H^0(\TT)=0$ we obtain that the square of the determinant line bundle
$\det R\Ga(\TT)$ is trivial. We have  natural isomorphisms 
(see \cite{Del})
$$\det R\Ga(\TT)^2\simeq \det R\Ga(\TT)\otimes \det R\Ga(\om\TT^{-1})\simeq
\det R\Ga(\om)\otimes\lan\TT,\om\TT^{-1}\ran^{-1},$$
where $\lan L,M\ran$ is Deligne's symbol.
We can compute the $4$-th power of this symbol:
\begin{align*}
\lan\TT,\om\TT^{-1}\ran^4&\simeq\lan\om(-p_1-p_2),\om(p_1+p_2)\ran\simeq
\lan\om,\om\ran\otimes\lan\OO(p_1+p_2),\OO(p_1+p_2)\ran^{-1}\\
&\simeq\lan\om,\om\ran\otimes\OO(-p_1)|_{p_1}\otimes\OO(-p_2)|_{p_2}\simeq
\lan\om,\om\ran\otimes\om|_{p_1}\otimes\om|_{p_2}.
\end{align*}
We conclude that
$$\det R\Ga(\TT)^8\simeq \det R\Ga(\om)^4\otimes\lan\om,\om\ran^{-1}\otimes
\om|_{p_1}\otimes\om|_{p_2}.$$
Thus, from the triviality of $\det R\Ga(\TT)^2$ on the complement to $Z$
we obtain the triviality of the canonical bundle on $S-Z$. Since
$\pi^{-1}(U)$ is an open subset of $S-Z$ which is \'etale over $U$
of degree $2^{2g}$, we see that the $2^{2g}$-th power of $\om_U$ is trivial. 
\ed 

\begin{cor} For a smooth projective curve $C$ of genus $g$ and a point
$p_0\in C$ let us consider the set 
$X_{p_0}\subset C-p_0$ consisting of points $p$ such that
there exists an effective divisor $D$ with $\OO_C(2D+p+p_0)\simeq\om_C$.
Then there exists a divisor $E$ supported on $X_{p_0}\cup\{ p_0\}$ 
such that $\OO_C(E)\simeq\om_C^{\otimes 2^{2g}}$.
\end{cor}

\begin{rem} The power $2^{2g}$ in the above corollary is probably not optimal.
For example, in the case $g=3$ it can be replaced by $8$.
\end{rem}


\begin{thebibliography}{99}
\bibitem{BFM} P. Baum, W. Fulton, R. MacPherson, 
{\em Riemann-Roch for singular varieties}, 
 Inst.\ Hautes Itudes Sci.\ Publ.\ Math.\ \textbf{45} (1975), 101--145. 

\bibitem{BM}
K.\,Behrend, Yu.\,Manin, \emph{Stacks of stable maps and
Gromov-Witten invariants,} Duke Math.\ J.\ \textbf{85} (1996), 
1--60.

\bibitem{Del} P.\,Deligne, {\em
Le determinant de la cohomologie\/},
Current trends in arithmetical algebraic geometry (Arcata, Calif., 1985),
Contemp. Math., {\bf 67} 93--177. 

\bibitem{Del2}P.\,Deligne, {\em Cohomologie \'a support propres},
SGA 4 XVII, Lecture Notes in Math. {\bf 305}, Springer, 1973.

\bibitem{F} W.\,Fulton, {\em Intersection theory}, Springer, 1998.

\bibitem{J2}  T.\,J.\,Jarvis, \emph{Torsion-free sheaves and moduli of
    generalized spin curves,} Compositio Math.\ \textbf{110} (1998),
    291-333.

\bibitem{J} 
T.\,J.\,Jarvis, \emph{Geometry of the moduli of higher spin
curves,} to appear in Internat.\ J.\ of Math.\ \textbf{11} (2000),
 \texttt{math.AG/9809138}.


\bibitem{JKV} T.\,J.\,Jarvis, T.\,Kimura, A.\,Vaintrob, \emph{Moduli spaces of
higher spin curves and integrable hierarchies}, to appear in Compositio
Math., {\tt math.AG/9905034}.

\bibitem{JKV2} T.\,J.\,Jarvis, T.\,Kimura, A.\,Vaintrob, \emph{Gravitational
descendants and the moduli space of higher spin curves},  
preprint,  {\tt math.AG/0009066}.

\bibitem{Ko}
M.\,Kontsevich, \emph{Intersection theory on the moduli space of curves
and the matrix Airy function,}
Commun.\ Math.\ Phys.\ \textbf{147} (1992), 1--23.

\bibitem{KM} 
M.\,Kontsevich, Yu.\,I.\,Manin, \emph{Gromov-Witten classes,
quantum cohomology,   and enumerative geometry}, Commun.\ Math.\ Phys.\
\textbf{164} (1994), 525--562. 

\bibitem{W3}
E.\,Witten, \emph{Two-dimensional gravity and intersection theory on moduli
space}, Surveys in Diff.\ Geom.\ \textbf{1} (1991), 243--310.

\bibitem{W1}
E.\,Witten, \emph{The $N$-matrix model and gauged WZW models,} Nucl.\ Phys.\ B
\textbf{371} (1992), no.\ 1-2, 191-245.

\bibitem{W2}
E.\,Witten, \emph{Algebraic geometry associated with matrix models of two
  dimensional gravity}, Topological methods in modern mathematics (Stony
Brook, NY, 1991), Publish or Perish, Houston, 1993, 235-269.


\end{thebibliography}
\end{document}